\documentclass{amsproc}
\usepackage[utf8]{inputenc}
  \usepackage{amsmath}
\usepackage{enumerate}
  \usepackage{amssymb}
  \usepackage{amsthm}
  \usepackage{bm}
  \usepackage{graphicx}
  \usepackage{color}
\usepackage{epsfig}
\usepackage{amsmath,amssymb}
\hoffset - 0.5 in
\theoremstyle{definition}

\theoremstyle{remark}

\numberwithin{equation}{section}

\newcommand{\tend}[3][]{\xrightarrow[#2\to#3]{#1}}

\newcommand{\egdef}{\stackrel{\textrm {def}}{=}}
\newcommand{\ds}{\displaystyle}

\title[Lebesgue component of multiplicity one]{A non-singular transformation whose spectrum has Lebesgue component of multiplicity one  }

\author{e. H. el Abdalaoui}
\address{Normandie University, University of Rouen
  Department of Mathematics, LMRS  UMR 60 85 CNRS\\
Avenue de l'Universit\'e, BP.12
76801 Saint Etienne du Rouvray - France .}
\email{elhoucein.elabdalaoui@univ-rouen.fr}
\urladdr{http://www.univ-rouen.fr/LMRS/Persopage/Elabdalaoui/}
\author{M. G. Nadkarni}
\address{Department of Mathematics, University of Mumbai, Vidyanagari, Kalina,  Mumbai, 400098, India}
\email{mgnadkarni@gmail.com}
\urladdr{http://insaindia.org/detail.php?id=N91-1080}
\dedicatory{}

\begin{document}
\maketitle
\begin{abstract} In this note we give an example of an ergodic non-singular map whose unitary operator admits a Lebesgue component of multiplicity one in its spectrum.
\end{abstract}

\section{Introduction}\label{intro}
  The purpose of this paper is to give an example of an ergodic non-singular transformation
  whose associated unitary operator admits Lebesgue component of multiplicity one  in its spectrum. M. Guenais \cite{Guenais} and Downarowicz-Lacroix \cite{Down} have shown that there exist measure preserving
  transformations with such property provided there  exists a sequence of analytic trigonometric polynomials $P_n, n=1,2,\cdots$ whose absolute values $\mid P_n\mid, n=1,2,\cdots$  converge to 1 in  some
  sense, and such that for each $n$, the coefficients of $P_n$ are real and equal in absolute value. In \cite{Guenais}, the initiating paper,  M. Guenais requires existence of such a sequence $P_n, n=1,2,\cdots$ whose $L^1$ norms converge to 1. Downarowicz and Lacroix requires the existence of the so called Barker polynomials of arbitrary high length. If we dispense with the requirement that the transformations be  measure preserving but require only that they be non-singular and ergodic then the problem has a non-speculative and concrete solution, for partial sums of power series expansions finite Blaschke products
  with real (non-vanishing) zeros provide us with ultraflat polynomials $P_n, n=1,2,\cdots$ with real coefficients, where, for any $n$ the coefficients of $P_n, n=1,2,\cdots$ are not equal in absolute value, but which can be used to construct non-singular ergodic transformations with the desired spectral property. We also mention that Guenais has constructed a concrete egodic measure preserving group action whose unitary group admits Haar component of multiplicity one in its spectrum. Our result shows that in the non-singular category one can get integer action with this property. We refer the reader to our papers \cite{Abd},\cite{Abd1} for a fuller discussion on connection between flat polynomials, $H^p$ theory and spectral questions in ergodic theory.\\

\section{Simple Blaschke Product and Ultraflat Sequence of Polynomials}
Let $S^1$ denote the circle group and $dz$ the normalized Lebesgue measure on $S^1$.
  A sequence $P_n(z), n=1,2,\cdots$ of analytic trigonometric polynomials of $L^2(S^1,dz)$
   norm 1 is said to be ultraflat if the sequence $| P_n(z)|, n=1,2,\cdots$ converges uniformly to the constant function $1$ as $n\rightarrow \infty$. J. E.  Littlewood \cite{Littlewood} asked if there exists
   an ultraflat sequence where, for each $n$, coefficients of $P_n$ are equal in absolute value, a question which J-P. Kahane \cite{Kahane} answered in the affirmative. The coefficients of
  polynomials in the ultraflat sequence constructed by Kahane are complex. It is an open question if there exist an ultraflat sequence of polynomials $P_n, n=1,2,\cdots$ where for each $n$, coefficients of $P_n$ are real and equal in absolute value.\\

  Our purpose is served by an ultraflat sequence $P_n, n=1,2,\cdots$, with coefficients of $P_n$  real and bounded away from 1 in absolute value. The sequence of partial sums of power series of a finite Blaschke product with real (non-vanishing) zeros provides us with such ultraflat sequence. We will use only the simplest of such Blaschke products.\\

  Consider the linear fractional map $\displaystyle \frac{z-\alpha}{1 - \alpha z}$, where $\alpha$ is real positive and less than one. This function maps the unit circle onto the unit circle and analytic in the closed disk of radius $\frac{1}{ \alpha} > 1$, so has a power series expansion valid on $S^1$.
  \begin{eqnarray*}
   (z-\alpha)\sum_{j=0}^\infty\alpha^{j}z^{j} &=& \sum_{j=0}^\infty\alpha^jz^{j+1} -\sum_{j=0}^\infty\alpha^{j+1}z^j\\
  &=& -\alpha + \sum_{j=1}^\infty \alpha^{j-1}(1-\alpha^2)z^j.
\end{eqnarray*}
 It can be checked that the sum of the squares of the coefficients of this power series is 1, as it should be, since  it is a function of absolute value 1 on $S^1$. The sequence of finite  sums
 $$R_n(z) = -\alpha + \sum_{j=1}^{n-1}\alpha^{j-1}(1 -\alpha^2)z^j, n= 1,2,\cdots$$

when divided by their $L^2(S^1,dz)$ norms,  is therefore ultraflat. Since $\|  R_n\|_2\rightarrow 1$ as $n\rightarrow \infty$, we see that for all $n \geq 2$, coefficients of $\frac{R_n}{\mid\mid R_n\mid\mid_2}$  in absolute value are less than some $\lambda <1.$

Let $\ds Q_m = \frac{R_m}{\mid\mid R_m\mid\mid_2}$ and choose a subsequence  $Q_{m_k}, k=1,2,\cdots$ of $Q_m, m=1,2,\cdots$ such that $\mid\mid Q_{m_k}(z)\mid -1\mid < \frac{1}{2^k}, k =1,2,\cdots$. Write $P_k = Q_{m_k}$. Clearly for any positive integer $l$,
$\mid \mid P_k(z^l)\mid - 1\mid < \frac{1}{2^k}$ and $P_k, k=1,2,\cdots$ are bounded away from zero on $S^1$. Let $h_0=1,h_j = m_1m_2\cdots m_j, j\geq 1$. We note that the finite products $\prod_{j=1}^k| P_j(z^{h_{j-1}})|,$  $k=1,2,\cdots$ converge uniformly to a non-vanishing function $g$
whose square  is also non-vanishing continuous function on $S^1$. We  claim that $\int_{S^1}g^2dz = 1$.
  Note that the constant term of $| P_j(z^{h_j})|^2$ is 1.
  Further the highest power of $z$ in  $P_1(z)P_2(z^{h_1})\cdots P_k(z^{h_{k-1}})$ is $h_k-1$, and the powers of $z$ in $P_{k+1}(z^{h_k})$ are multiples of $h_k$ so the constant term of
$$\prod_{j=1}^k\mid P_j(z^{h_{j-1}})\mid^2 ~~~~~\eqno (1)$$ is also 1. So these products integrate to 1. Clearly, since products (1) converges to $g^2$ uniformly, $\ds \int_{S^1}g^2dz =1$.\\

\section{ Unitary operators $U_T$ and $V_\phi$}

We will now construct a non-singular rank one non-dissipative, ergodic transformation $T$
   on the unit interval (equipped with Lebesgue measure $\nu$) and a function $\phi$ on $[0,1]$ taking values  $-1$ and $1$ such that the maximal spectral type of the unitary operator $V = V_\phi$ defined by
   $$(V_\phi f)(x)   = \phi (x)\Big(\frac{d\nu \circ T}{d\nu}(x)\Big)^{\frac12} f(Tx), f \in L^2([0,1], \nu)$$
has simple spectrum with maximal spectral type $g^2dz$. The transformation $T$ will be a generalized odometer.

\section{ Ultraflat Polynomials $P_k, k=1,2,\cdots$ and rank one transformation}

Let
$$P_k(z) = \sum_{j=0}^{m_k-1} a_{j,k}z^j.$$  Since $\mid\mid P_k\mid\mid_2^2 =1,$
$\ds \sum_{j=0}^{m_k-1} a_{j,k}^2 = 1.$  Write $p_{j,k} = a_{j,k}^2$.

We will now a construct a rank one non-singular transformation of generalized odometer type,
equivalently, a rank one transformation where no spacers added at any stage of  construction.
 This is done by the method of cutting and stacking \cite{Friedman} as follows. Let $\Omega_0 = \Omega_{0,0}$ denote the unit interval. At stage one of the construction we divide $\Omega_0$ into $m_1$ pairwise disjoint intervals, $\Omega_{0,1}, \Omega_{1, 1} \cdots, \Omega_{m_1-1, 1}$, of lengths $p_{0,1}, p_{1,1},\cdots, p_{m_1-1, 1}$ respectively.  For each $j, 0\leq j \leq m_1-2$, we map the interval $\Omega_{j,1}$ linearly onto $\Omega_{j+1,1}$. We view the intervals $\Omega_{0,1}, \Omega_{1,1}\cdots, \Omega_{m_1-1, 1}$ as stacked one above the other.
We thus get a stack of certain height $h_1 = m_1$, together with a map $T$ which is defined on all intervals of the stack except the interval $\Omega_{m_1-1,1}$ at the top of the stack.  Note that since $p_{j,1} \neq p_{j+1,1}$, $T$ will not be measure preserving. This completes the first stage of the construction.\\

Suppose we have obtained at the end of $(k-1)$th stage of construction a stack of height
$h_{k-1} = m_1m_2\cdots m_{k-1}$ and a map $T$ which maps each interval of the stack
 linearly onto one above it, except that $T$ remains undefined on the top piece of the stack. Let $\Omega_{0,k-1}$ denote the bottom interval of this stack. $T^{h_{k-1}-1}\Omega_{0,k-1}$ is then the top of the stack.
At $k$-th stage of construction we divide $\Omega_{0, k-1}$ into $m_k$ intervals $\Omega_{0,k},\Omega_{1,k}, \cdots, \Omega_{m_k-1, k}$ in the ratios $p_{0,k}, p_{1,k}, \cdots, p_{m_k-1, k}$, and map, for $0 \leq j \leq m_k-2$, the interval $T^{h_{k-1}-1}\Omega_{j,k}$ linearly onto $\Omega_{j+1,k}$. We thus get, at the end of $k$-th stage of construction a stack of height
$h_k = m_1m_2\cdots m_k$, where each piece of the stack is mapped  by $T$ linearly onto the one above it except that $T$ remains undefined on the top piece of the stack which is $T^{h_{k-1}-1}\Omega_{m_k-1,k}$. Since $p_{i,k}$'s are  bounded away from 1, we see that measure of
$T^{h_{k-1}-1}\Omega_{m_k-1,k}$ goes to 0 as $k$ tends to $\infty$, whence  $T$ is eventually defined almost everywhere and we get a non-singular transformation $T$ on $[0,1]$ which we call non-singular generalized odometer. Since $p_{k, j}$'s are bounded away from 0 over all $k,j$, the maximum length of the intervals of $k$-th stack goes to 0 as $k$ tends to $\infty$. So a Lebesgue density argument allows us to prove that $T$ is ergodic. We omit this proof.\\

 We also remark that there is no finite measure on Borel subsets of $[0,1]$ which is invariant under $T$ and mutually absolutely continuous with respect to $\nu$. This in turn implies that the unitary operators $U_T$ and $V_\phi$ defined in the next section do not admit eigenvalues.\\

The action $f \rightarrow f\circ T$ does admit a countable dense subgroup of $S^1$ as eigenvalues with measurable eigenfunctions of absolute value 1 and which separate points
$[ 0,1]$ module $\nu$-null sets. Let $\Gamma$ denote this group of eigenvalues of $T$. For each $\gamma \in \Gamma$  we can choose an eigenfunction $e_\gamma$ which is measurable, of absolute value 1, and satisfies  for all $\gamma_, \gamma_2 \in \Gamma$, $e_{\gamma_1\gamma_2} = e_{\gamma_1}e_{\gamma_2}$ modulo $\nu$ null sets. We will use these facts in the next section to prove that the maximal spectral type of $U_T$ is quasi-invariant and ergodic under translation action of $\Gamma$ on $S^1$.\\

\section{ $V_\phi$,  Lebesgue nature of its Spectrum}

Let $\phi$ denote a function on $[0,1]$, constant on each level of any stack except the top level, assuming values $-1$ and $1$. This function will be  constructed inductively as we proceed.

Let $T$, $\phi$ be as above. On $L^2([0,1],\nu)$ define
$$(U_T f)(x)   \egdef \Big(\frac{d\nu\circ T}{d\nu}(x)\Big)^{1/2} f(Tx), f \in L^2([0,1], \nu),$$
$$(V_\phi f)(x)  \egdef (Vf)(x) = \phi(x)\cdot (U_Tf)(x), f \in L^2([0,1], \nu).$$
$U_T$, and $V$ are unitary operators. Since one has the formula

$$\frac{d\nu\circ T^2}{d\nu} = \frac{d\nu\circ T}{d\nu}\frac{d\nu\circ T}{d\nu}\circ T,$$ and more generally for $n \geq 0$

$$\frac{d\nu\circ T^n}{d\nu} = \frac{d\nu\circ T}{d\nu}\frac{d\nu\circ T^{n-1}}{d\nu}\circ T.$$
We have for $n \geq 0$,
$$(U^n_Tf)(x) = \Big(\frac{d\nu\circ T^n}{d\nu}(x)\Big)^{1/2}f(T^nx),$$

$$(V^nf)(x)=\Big(\prod_{j=0}^{n-1}\phi(T^j(x))\Big)\Big(\frac{d\nu\circ T^n}{d\nu}(x)\Big)^{1/2}f(T^nx),$$

Since any bounded measurable function can be approximated in $L^2([0,1], \nu)$ by finite linear combinations of $V^j1_{\Omega_{0,n}}, 0 \leq j \leq h_n-2$ for large enough $n$, $V$ has multiplicity one.\\

Write $Tf = f\circ T$. Using facts such as $(AB)^{-1} = B^{-1}A^{-1}$ for invertible operators $A, B$, and that multiplication operators commute,  we have from the formula for $V^nf$ given above:
$$(V^{-n}f)(\cdot) = T^{-n}\circ\Big(\Big(\prod_{j=0}^{n-1}\phi(T^j(\cdot)\Big)\Big)^{-1}\Big(\frac{d(\nu\circ T^n)}{d\nu}(\cdot)\Big)^{-1/2}f(\cdot))$$
$$=\Big(\prod_{j=0}^{n-1}\phi (T^{j-n}(\cdot))\Big)^{-1}\Big(\frac{d\nu\circ T^n}{d\nu}(T^{-n}(\cdot))\Big)^{-1/2}f(T^{-n}(\cdot)),$$
whence
$$(T^{-n}f)(\cdot) =  \prod_{j=0}^{n-1}\phi (T^{j-n}(\cdot))\Big(\frac{d(\nu\circ T^n)}{d\nu}\Big)^{1/2}(T^{-n}(\cdot))(V^{-n}f)(\cdot).$$

\noindent{}  Recall that $\Omega_{0,k-1}$ is the base of the stack of height $h_{k-1}$ and $\Omega_{0,k}$, $\Omega_{1,k}, \cdots, \Omega_{m_k-1, k}$ is a partition of $\Omega_{0,k-1}$. We have
$$\Omega_{0,k-1} = \bigcup_{j=0}^{m_k-1}T^{R_{j,k}}(\Omega_{0,k}),$$
where $R_{j,k} = jh_{k-1}$ is the $j$-th return time of a point in $\Omega_{0,k}$ into $\Omega_{0,k-1}$.
\begin{eqnarray*}
1_{\Omega_{0,k-1}} &=& \sum_{j=0}^{m_k-1}1_{\Omega_{0,k}}\circ T^{-R_{j,k}}\\
&=& 1_{\Omega_{0,k}} + \sum_{j=1}^{m_k-1}c_{j,k}\Big(\frac{d\nu\circ
T^{R_{j,k}}}{d\nu}\circ T^{-R_{j,k}}\Big)^{1/2}(\cdot)(V^{-R_{j,k}}1_{\Omega_{0,k}})(\cdot)\\
&=& W_k(V)1_{\Omega_{0,k}},
\end{eqnarray*}
where
$$W_k(z) = 1 + \sum_{j=1}^{{m_k-1}}c_{j,k}\Big(\frac{d\nu\circ T^{R_{j,k}}}{d\nu}(T^{-R_{j,k}}))\Big)^{1/2}(\cdot)(z^{-R_{j,k}}),$$

 \noindent{} and $c_{j,k} =\ds \prod_{i=0}^{R_{j,k}-1}\phi (T^{i-R_{j,k}}(x)), x \in \Omega_{j,k}, 1 \leq j \leq m_k-1$, a constant of absolute value one. Note that the constants $c_{j,k}$'s can be preassigned and $\phi$ can be so defined that the above relation holds for all $(j,k)$.
We will define  $\phi$ inductively in such a way that $c_{j,k} = -1$, for all $j, k, 1\leq j \leq m_{k}-1$. This will ensure, after we have substituted  values of the Radon-Nikodym derivatives, that
$$W_k(z) = -\frac{1}{\alpha}P_k(z^{h_{k-1}}).$$

\noindent{} Choose $\phi = -1$ on $\Omega_{0,1}$ and equal to
$1$ on $\Omega_{j,1}, 1\leq j \leq m_1-2$, and verify  that
$$c_{j,1} = \prod_{i=0}^{R_{j,1}-1}\phi(T^{i-R_{j,1}}x) = -1,~~~ x \in \Omega_{j,1},~~ 1\leq j \leq m_1 -1.$$
Assume now that $\phi$ has been defined on first $h_{k-1}-1$ levels of the stack of height
$h_{k-1}$, $\phi$ is constant on each level with value $-1$ or $1$, such that for $l\leq k-1$
$$c_{j,l} = -1,~~1\leq j \leq m_l-1.$$
Let $b = \prod_{i=0}^{h_{{k-1}}-2}\phi (T^i(x)), x\in \Omega_{0, k-1}$, a constant independent of $x\in \Omega_{0,k-1}$, of value $-1$ or $1$. We now define $\phi$ on the intervals
$T^{h_{k-1}-1}\Omega_{j,k}, 0 \leq j \leq m_k-2$ as follows: if $b=1$ the value of $\phi$ is  -1 on $T^{h_{k-1}-1}\Omega_{0,k}$, and 1 on the intervals $T^{h_{k-1}-1}\Omega_{j,k}, 1\leq j\leq m_k-2$.
if $b=-1$ the value of $\phi$ is  1 on $T^{h_{k-1}-1}\Omega_{0,k}$, and -1 on the intervals $T^{h_{k-1}-1}\Omega_{j,k}, 1\leq j\leq m_k-2$.
This ensures that $c_{j,k} = -1,$ $1\leq j \leq m_k-2$. Thus we have defined $\phi$ inductively on all of $[0,1]$ so that $c_{j,k} =-1,$  $1\leq j \leq m_{k}-1, k = 1,2 \cdots$.

We now observe that for $x\notin T^{R_{j,k}}\Omega_{0,k}$,
$$V^{-R_{j,k}}1_{\Omega_{0,k}}(x) = 0,$$  and that for $x\in T^{R_{j,k}}\Omega_{0,k}$,
$$\frac{d\nu\circ T^{R_{j,k}}}{d\nu}(T^{-R_{j,k}}(x)) = \frac{p_{j,k}}{p_{0,k}}.$$
We thus have, with $c_{0,k} =1, R_{0,k} =0$,
$$1_{\Omega_{0,k-1}}  =\sum_{j=0}^{{m_k-1}}c_{j,k}\Big(\frac{p_{j,k}}{p_{0,k}}\Big)^{1/2}(V^{-R_{j,k}}1_{\Omega_{0,k}})(\cdot).$$
Let us normalize $1_{\Omega_{0,k}}$ and write
$$f_k =\Big (\frac{1}{ \nu(\Omega_{0,k})}\Big)^{1/2}1_{\Omega_{0,k}} = \Big( \frac{1}{\prod_{j=1}^{k} p_{0,j}}\Big)^{1/2}1_{\Omega_{0,k}},$$
$$f_{k-1} = \big(p_{0,k}\big)^{1/2}\Big(1 + c_{1,k}\Big(\frac{p_{1,k}}{p_{0,k}}\Big)^{1/2}V^{-R_{1,k}} + \cdots + c_{m_k-1,k}\Big(\frac{p_{m_k-1,k}}{p_{0,k}}\Big)^{1/2}V^{-R_{m_k-1,k}}\Big)f_k.$$
Now $\nu(\Omega_{0,0}) = 1$ so $f_0 = 1_{\Omega_{0,0}}$. We have by iteration
$$f_0 = \Big(\prod_{j=1}^kS_j(V)\Big)f_k,$$
where
$$S_j(z) = \Big(p_{0,j}\Big)^{1/2}\Big(1 + c_{1,j}\Big(\frac{p_{1,j}}{p_{0,j}}\Big)^{1/2}z^{-R_{1,j}} + \cdots + c_{m_j-1,j}\Big(\frac{p_{m_j-1,j}}{p_{0,j}}\Big)^{1/2}z^{-R_{m_j-1,j}}\Big).$$
Note that by choice of $p_{i,j}$ and $c_{i,j}$,
$$S_j(z)  = -P_j(z^{h_{j-1}}), j=1,2,\cdots, \mid S_j(z)\mid^2 = \mid P_j(z^{h_{j-1}})\mid^2.$$

Let $V^n = \int_{S^1}z^{-n}dE, n \in \mathbb Z$, be the spectral resolution of the unitary group $V^n, n \in \mathbb Z$, and

$$ (V^nf_k, f_k) = \int_{S^1}z^{-n}(E(dz)f_k,f_k) =\int_{S^1}z^{-n}d\sigma_k,$$

where $\sigma_k (\cdot) = (E(\cdot)f_k, f_k)$; the maximal spectral type of $E$ is given by $\vee_{k=0}^\infty \sigma_k$.\\

 Since $T^i\Omega_{0,k}, 0 \leq i\leq h_k-1$\
are pairwise disjoint, we see that the sequence $\sigma_k, k=1,2,\cdots $ converges weakly to $dz$ as $k\rightarrow \infty$.

We have for all integers $l$

$$(V^lf_0,f_0) = \int_{S^1}z^{-l}d\sigma_0 =\int_{S^1}z^{-l} \prod_{j=0}^k\mid P_j(z^{h_{j-1}})\mid^2d\sigma_k,$$

whence
$$d\sigma_{0} = \prod_{j=1}^k\mid P_j(z^{h_{j-1}})\mid^2d\sigma_k.$$

Since the spectral measure $E$ is atomfree and $P_j, 1 \leq j \leq k$, can vanish only at finitely points
we see that for all $k$, $\sigma_0$ and $\sigma_k$ are mutually absolutely continuous, whence $\sigma_0$ is the maximal spectral type of $V_\phi$.\\

Now $\sigma_k \rightarrow dz$ weakly as $k\rightarrow \infty$ and $\prod_{j=1}^k|P_j(z^{h_{j-1}})|^2 \rightarrow g^2$ uniformly, whence $ \prod_{j=1}^k| P_j(z^{h_{j-1}})|^2 d\sigma_k
\rightarrow g^2dz$ weakly as $k\rightarrow \infty$. Thus $\sigma_0$ is the measure $g^2dz$.
Thus the maximal spectral type of $V_\phi$ is mutually absolutely continuous with respect to the  Lebesgue measure.\\

We will now show that the maximal spectral type of $U_T$ is quasi-invariant and ergodic
with respect to the translation action by $\Gamma$ on $S^1$. Consider the unitary group
$J_\gamma: f \rightarrow e_\gamma f, f \in L^2([0,1], \nu), \gamma \in \Gamma$, which satisfies with $U_T$ the Weyl commutation relation $U_TJ_\gamma = \gamma J_\gamma U_T$.
Since the functions $e_\gamma, \gamma \in \Gamma$, separate points and $T$ is ergodic with respect to $\nu$
the pair of unitary groups $(U_T^n, n\in \mathbb Z, J_\gamma, \gamma \in \Gamma)$, is irreducible in the sense that only closed subspaces of $L^2([0,1], \nu)$ invariant under all $U_T^n, n \in\mathbb Z$, and all $J_\gamma, \gamma \in \Gamma$, are the trivial ones. This in turn implies that the spectral measure $F$ of $U_T$ satisfies with $J_\gamma, \gamma \in \Gamma$, the relation
$J_\gamma F(\cdot) J_{\gamma^{-1}} = F((\cdot)\gamma), \gamma \in \Gamma$, and  the system $(F, J_\gamma, \gamma \in \Gamma)$ is irreducible. It follows from this that the maximal spectral type of $F$, say $\sigma$, is ergodic and quasi-invariant under $\Gamma$. Since $\Gamma$ is dense in $S^1$ and $\sigma$ is ergodic under $\Gamma$ we see that $\sigma$ is either equivalent to $dz$ on $S^1$ or singular to it. (\cite{Nad},11.11, 12.14,13.3, 13.4, 13.5).\\

 Next we show that the maximal spectral type of $U_T$ is singular to $(dz)$. Recall  that the constant term of $a_{0,j}$ of $P_j(z)$ is negative and converges to $-\alpha$ as $j\rightarrow \infty$ and that $\ds P_j(z) \rightarrow \frac{z-\alpha}{1 - \alpha z}$ as $j\rightarrow \infty$. The polynomial
$L_j(z)  \egdef 2| a_{0,j}|+ P_j(z)$ have $L_2(S^1,dz)$ norm 1 and
$$ \ds | L_j(z)| \rightarrow \Big| \frac{z-\alpha}{1 - \alpha z} + 2\alpha\Big|$$ uniformly as $j \rightarrow \infty$. The function $\ds \Big| \frac{z-\alpha}{1 - \alpha z} + 2\alpha\Big|$ is non-constant, it has $L^2(S^1,dz)$ norm 1 since each $L_j$ has $L^2(S^1, dz)$ norm 1. So its $L^1(S^1,dz)$, norm, say $l$, is strictly less than one, whence, as $j\rightarrow \infty$, $\ds \int_{S^1}| L_j(z)| dz \rightarrow l < 1$.   \\

  {\bf Lemma.} {\it If the maximal spectral type of $U_T$  has the same null sets as $dz$ then $\ds \int_{S^1} |L_n(z)| dz\rightarrow 1 $  as $n\rightarrow \infty$.}\\

(Note that this contradicts what is stated just above, so the maximal spectral type of $U_T$ is singular to $dz$ by the dichotomy result.)\\

\begin{proof} The maximal spectral type of $U_T$ is given by the weak limit of measures
$\prod_{j=1}^n| L_j(z^{h_{j-1}})|^2dz, n=1,2,\cdots.$ It is known that   $\prod_{j=1}^n| L_j(z^{h_{j-1}})|, n=1,2,\cdots$ converge over a subsequence, say $n_k, k=1,2,\cdots$, to $\ds \sqrt\frac{d\sigma}{dz}$ a.e. $(dz)$, (see \cite{Abd}, Proposition 5.2).
Assume that $\sigma$ and $dz$ have the same null sets. Then there is a subsequence $n_k, k=1,2,\cdots$ such that, for almost all $z$, we have
$$\prod_{j=1}^{n_k}|L_j(z^{h_{j-1}})|(z)\tend{k}{+\infty}\sqrt{\frac{d\sigma}{dz}}(z).$$
Since $\ds \frac{d\sigma}{dz}(z)$ is positive a.e. $(dz)$ the partial products
 $\ds \prod_{j=n_k+1}^{n_{k+1}}|L_j(z^{h_{j-1}}|, k=1,2,\cdots $ converge almost everywhere $(dz)$ to $1$.
We thus get,
$$\int_{S^1} \prod_{j=n_k+1}^{n_{k+1}}|L_j(z^{h_{j-1}})| dz \tend{k}{+\infty}1,$$
since, by Cauchy-Schwarz inequality, for all $k$,
$$\int_{S^1}\prod_{j=n_k+1}^{n_{k+1}}|L_{j}(z^{h_{j-1}})|dz \leq 1.$$

\noindent{}Let $M_k$ be any $L_j, n_{k}+1\leq j \leq n_{k+1}$. We will show that
$$\int_{S^1}\prod_{j=n_k+1}^{n_{k+1}}|L_{j}(z^{h_{j-1}})|dz \leq \Big(\int_{S^1}\mid M_k(z)\mid dz\Big)^{\frac{1}{2}},$$
wherein letting $k \rightarrow \infty$ the Lemma follows.
To fix ideas we take $M_k = L_{n_{k+1}}$.
Now, by Cauchy-Schwarz inequality, we have

$$\int_{S^1} \prod_{j=n_k+1}^{n_{k+1}}|L_j(z^{h_{j-1}})| dz = \int_{S^1} \sqrt{|L_{n_{k+1}}(z^{h_{n_{k+1}-1}})|} \sqrt{|L_{n_{k+1}}(z^{h_{n_{k+1}-1}})|} \prod_{j=n_k+1}^{n_{k+1}-1}|L_j(z^{h_{j-1}})| dz$$
$$\leq \Big(\int_{S^1} |L_{n_{k+1}}(z^{h_{n_{k+1}-1}})| dz \Big)^{\frac12}
\Big(\int_{S^1} |L_{n_{k+1}}(z^{h_{n_{k+1}-1}})| \prod_{n_k+1}^{n_{k+1}-1}|L_j(z^{h_{j-1}})|^2 dz\Big)^{\frac12}$$
$$\leq \Big(\int_{S^1} |L_{n_{k+1}}(z^{h_{n_{k+1}-1}})| dz \Big)^{\frac12},$$

since, again by Cauchy-Schwarz inequality, we have

$$\int_{S^1} |L_{n_{k+1}}(z^{h_{n_{k+1}-1}})| \prod_{j=n_k+1}^{n_{k+1}-1}|L_j(z^{h_{j-1}})|^2 dz$$
$$= \int_{S^1} \prod_{j=n_k+1}^{n_{k+1}}|L_j(z^{h_{j-1}})|  \prod_{j=n_k+1}^{n_{k+1}-1}|L_j(z^{h_{j-1}})| dz$$
$$\leq \Big(\int_{S^1} \prod_{j=n_k+1}^{n_{k+1}}|P_j(z^{h_{j-1}})|^2 dz\Big)^{\frac12}
\Big(\int_{S^1} \prod_{j=n_k+1}^{n_{k+1}-1}|L_j(z^{h_{j-1}})|^2 dz\Big)^{\frac12}=1.$$

The lemma follows.
\end{proof}

We now construct the desired non-singular non-dissipative ergodic transformation $\tau$ with
Lebesgue component in the spectrum of $U_\tau$. Since $V_\phi$ has spectrum different from that of $U_T$, $\phi$  is not a multiplicative coboundary, i.e, it is not of the form $\frac{\xi\circ T}{\xi}$, for any measurable function $\xi$. Consider now the space $X = [0,1]\times \{-1,1\}$, equipped with product measure, denoted by $\mu$, where $[0,1]$ has Lebesgue measure and $\{-1,1\}$ has uniform probability distribution. View $\{-1,1\}$ as a multiplicative group with two character $\xi_0, \xi_1$, with $\xi_0$ the identity character. Define $\tau (x,y) = (Tx, \phi (x)y), (x,y) \in X$.
The transformation $\tau$ is non-singular with respect to $\mu$, and ergodic since $\phi$ is not a multiplicative coboundary. Define
$$U_\tau f = \sqrt{\frac{d\mu\circ \tau}{d\mu}}f\circ \tau, f \in L^2(X,\mu).$$
The space $L^2(X,\mu)$ is orthogonal sum of $L^2([0,1],\nu)$ and $\xi_1L^2([0,1], \nu)$
each invariant under $U_\tau$. The restriction of $U_\tau$ to $L^2([0,1], \nu)$ is $U_T$
 and its  restriction to $\xi_1L^2([0,1],\nu)$ is $V_\phi$ since  $\xi_1 (\phi) = \phi$.
Clearly $U_\tau$ admits Lebesgue component of multiplicity one.\\

{\bf Remarks}
 \begin{enumerate}[1)]
\item There is no ultraflat sequence of polynomials with non-negative coefficients which are bounded away from one. For, if $\ds P_n(z) = \sum_{j=0}^{m_n}a_{j, n}z^j, n=1,2,\cdots,$ is an ultraflat sequence of polynomials such that for each $n$, $0\leq a_{j,n} < \lambda <1$, $ 0 \leq j \leq m_n$ then since  it is ultraflat $$\sum_{j=0}^{m_n}a_{j,n}^2 = 1
<  \sum_{j=1}^{m_n}a_{j,n} = P_n(1),$$ whence
$$P_n(1) - 1 = \sum_{j=0}^{m_n}(a_{j,n}-a_{j,n}^2) > 1 -\lambda, $$
contradicting that $P_n(1) \rightarrow 1$ as $n\rightarrow \infty$.

  \item We do not know if there exists a sequence $P_n, n=1,2,\cdots$ of analytic polynomials of $L^2(S^1,dz)$ norm 1, with coefficients non-negative and bounded away from 1, such that $\mid P_n(z)\mid \rightarrow 1$ a.e $(dz)$ as $n\rightarrow \infty$. An affirmative answer will imply that there are ergodic, non-dissipative, non-singular transformations $\tau$ with simple  Lebesgue spectrum for $U_\tau$, answering a question of Banach in the non-singular category \cite{Abd}, \cite{Abd1}, \cite{CN}. In particular it is not known if for a non-dissipative non-singular generalized odometer action $T$, $U_T$ always has singular spectrum.  Using the Peyri\`ere-Brown method \cite{pey},\cite{Br}, it is can be shown that the spectrum of the generalized odometer is singular if the following holds
$$\sum_{k \geq 1} \frac{\big| \sum_{i \neq j} \sqrt{p_{i,k}p_{j,k}}\big|^2}{2(m_k-1)}=+\infty.$$
 \item The Banach problem asks if there is a Lebesgue measure preserving transformation on $\mathbb R$ which has simple Lebesgue spectrum. This problem is stated in 1959's Ulam book \cite[p.76]{Ulam}. A similar problem is mentioned by Rokhlin in \cite{Rokh}. Precisely, Rokhlin asked on the existence of an ergodic measure preserving transformation on a finite measure space whose spectrum is Lebesgue type with finite multiplicity. Later, in 1966, Kirillov in \cite{Kiri} wrote ``there are grounds for thinking that such examples do not exist". However he has described a measure preserving  action (due to M. Novodvorskii) of the group $(\bigoplus_{j=1}^\infty{\mathbb {Z}})\times\{-1,1\}$ on the compact dual of discrete rationals whose unitary group has Haar spectrum of multiplicity 2. Similar group actions with higher finite even multiplicities are also given.

\item Finite measure preserving  $T$ with $U_T$ having Lebesgue component of finite even multiplicity have been constructed by J. Mathew and M. G. Nadkarni\cite{MN}, M. Queffelec\cite{Q}, and O. Ageev\cite{Ag}.

\item For the case of the non-singular maps, Ismagilov established the connection between the classical Riesz product and the spectral type of some Mackey actions \cite{Ism1}, \cite{Ism2}, \cite{Ism3}. For all those examples, he proved  that the spectrum is singular.
\end{enumerate}

  \end{document}